\documentclass[twoside,11pt]{article}
\usepackage{jmm}
\usepackage{booktabs}
\nolinenumbers
\begin{document}
\Title{Optimal Control of a Bioeconomic Crop-Energy System with Energy Reinvestment}

\Author{Othman Cherkaoui Dekkaki$^\dag$}

\Address{$^{\dag}$ College of Computing, University Mohammed VI
Polytechnic, Benguerir, Morocco.}

\Email{othman.cherkaoui-dekkaki-ext@um6p.ma}
\Markboth{O. Cherkaoui Dekkaki}{Optimal Control of a Bioeconomic Crop-Energy System with Energy Reinvestment}

\Abstract{We develop a continuous-time optimal control model for allocating agricultural crop residues between bioenergy production and soil fertility restoration. The system includes a circular reinvestment channel: a portion of the accumulated bioenergy stock is reinvested to enhance soil fertility, thereby closing the loop between energy use and ecological regeneration. The dynamics are governed by a three-state system with a single allocation control. The objective is to maximize a discounted net benefit that accounts for energy revenue, soil value, and operational costs. We apply the Pontryagin Maximum Principle in current-value form to derive necessary optimality conditions and a bang–interior–bang control; no singular arc under our parameterization. Direct optimization confirms a quasi-turnpike phase and shows how the planning horizon and reinvestment efficiency shift switching times and interior duration. The results highlight the strategic role of energy reinvestment in sustainable residue management.
}


\section{Introduction}

Sustainable use of crop residues is at the intersection of environmental stewardship, agricultural productivity, and the development of renewable energy. As by-products of harvesting operations, residues such as straw, stalks, and husks represent a globally abundant biomass resource estimated at about $3.8$ billion tonnes per year (dry matter) \emph{worldwide} in the early 2000s \cite{lal2005}; see also \cite{smith2014afolu} for broader AFOLU context.
 Their potential applications range from soil fertility enhancement through organic matter recycling to energy generation via biochemical or thermal processes. However, these alternative uses are inherently competitive: returning residues to the soil supports long-term ecosystem health, while diverting them to energy production yields short-term economic gains. This trade-off is particularly critical in systems where soil organic carbon stocks are vulnerable and organic inputs are scarce.

To address such allocation challenges, mathematical modeling, particularly in the field of bioeconomics, has played a central role. The foundational work of of Clark \cite{clark2010mathematical}, Gordon \cite{gordon1954economic}, and others introduced dynamic frameworks that link biological processes to economic optimization. Subsequent models have extended these ideas to agricultural systems, incorporating soil dynamics, crop productivity, and market behavior \cite{flichman2011bioeconomic, antle2001econometric, van1994dynamic}. In parallel, policy-oriented assessments such as the US billion-ton update \cite{perlack2011us} and EU biomass reports have evaluated sustainable removal thresholds for crop residues to preserve soil functions. The circular bioeconomy literature further advances this view by promoting residue reintegration strategies such as composting, biochar application, or digestate recycling to close carbon and nutrient loops \cite{lehmann2015biochar, anikwe2023circular}.

Despite these developments, most assessments rely on static criteria or scenario analysis, with limited integration of dynamic feedback between ecological and economic components. This is where control theory provides a powerful framework. Optimal control methods, especially Pontryagin’s Maximum Principle (PMP) \cite{pont64, kamien2012dynamic, aseev7}, enable the formal derivation of time-dependent allocation strategies. In agricultural contexts, PMP has been applied to irrigation scheduling \cite{liu2016optimal}, pest control \cite{silva2017optimal}, and soil nutrient management \cite{yirga2010social}, but few models integrate biomass use, energy recovery, and soil regeneration into a unified optimization framework.

Recent contributions to bioeconomic modeling highlight the potential of PMP and related methods in capturing complex sustainability trade-offs. Studies have examined optimal resource allocation in ecological and technological systems, including microbial metabolite production \cite{caillau22tur}, species selection \cite{dj21turn, dj22op}, and anaerobic digestion \cite{rapaport19}. In the domain of waste-to-energy, previous work \cite{cherkaoui2022bio, cherkaoui2023acc} employed PMP to optimize investment and valorization strategies under environmental constraints, while \cite{cherkaoui2024ijdc} demonstrated how direct optimization methods can numerically resolve state-constrained control problems with feedback effects in waste recovery chains.

Building on this foundation, we propose a continuous-time optimal control model that allocates agricultural residues between energy production and soil fertility restoration. The model introduces a reinvestment mechanism whereby a fraction of the accumulated energy output enhances future soil productivity representing conservation spending, nutrient cycling, or infrastructure co-benefits. This circular feedback is embedded within a nonlinear three-state dynamical system and studied through the PMP framework. We derive necessary optimality conditions, characterize the structure of optimal trajectories, and explore the resulting strategies through direct numerical optimization, analyzing their sensitivity to reinvestment efficiency and planning horizon.

The contribution of this work is threefold:
\begin{itemize}
\item It advances control-based bioeconomic modeling by endogenizing ecological feedbacks from energy production to soil systems.
\item It operationalizes circular bioeconomy principles within a rigorous dynamic optimization framework, moving beyond static allocation rules.
\item It provides a generalizable methodology for evaluating sustainable residue management strategies under competing objectives.
\end{itemize}

\paragraph*{Related work and positioning.}
This paper contributes to the optimal control of residue–energy systems by introducing an explicit \emph{energy-to-soil reinvestment} state channel: a fraction of accumulated bioenergy is fed back to the soil, enhancing soil organic carbon (SOC) and, through it, crop productivity. Prior optimization studies in the biomass/waste domain largely treat residues as a supply and logistics problem aimed at cost-effective energy delivery or storage management, without a dynamic soil-feedback state that closes the loop from energy to fertility \cite{Nunes2024,sandoval2024}. Empirically, the loop we model has two well-documented sides: removing residues reduces SOC on average (a recent global review reports an $\approx 11\%$ decrease in topsoil SOC under residue harvest), while stabilizing carbon from residues (e.g., via biochar) tends to increase SOC substantially (meta-analysis mean $\approx 26\%$) \cite{alvarez2024,tian2024}. Embedding this evidence into a tractable three-state OCP allows us to study how the economic environment (energy and soil prices) and the reinvestment efficiency jointly shape the timing of residue diversion vs.\ soil return. Our finite-horizon results complement the broader literature on turnpike behavior where long-horizon optima spend most time near a steady state by showing when the reinvestment channel still yields bang–interior–bang patterns typical of systems with strong economic steady states \cite{trelatzuazua2025,geshkovski2022}.

This paper is structured as follows. Section~\ref{sec:model} introduces the mathematical formulation of the system, defining the state dynamics, control variable, and key bioeconomic assumptions. Section~\ref{sec:ocp} formally states the optimal control problem (OCP), including the objective functional and admissible trajectories. Section~\ref{sec:PMP} presents the analytical framework based on the Pontryagin Maximum Principle, characterizing the necessary optimality conditions and the switching behavior of the control. Section~\ref{sec:numerics} illustrates the numerical resolution of the OCP using direct optimization methods, and compares control policies under different reinvestment scenarios and planning horizons.

\section{Mathematical Formulation}
\label{sec:model}

We propose a continuous-time dynamical model that captures the allocation of agricultural crop residues between two competing pathways: the restoration of soil fertility and the generation of bioenergy. This formulation is motivated by a central trade-off in circular bioeconomy systems: diverting biomass to energy delivers immediate economic benefits but may undermine soil health, while returning biomass to the soil supports long-term productivity at the expense of forgone energy.

\textbf{A distinguishing feature of the model} is the inclusion of a feedback mechanism whereby cumulative bioenergy output contributes positively to future soil fertility. This reinvestment captures the indirect benefits of energy use such as infrastructure development, adoption of improved farming practices, or income-driven conservation and embeds a memory effect in the system dynamics. From a control-theoretic standpoint, this structure allows for the emergence of path-dependent behaviors and long-term policy impacts, making it well suited for sustainability analysis.

To maintain tractability and analytical transparency, we adopt a \textbf{linear-affine structure} in both state and control variables. This choice enables the application of Pontryagin's Maximum Principle and facilitates the derivation of qualitative properties of optimal solutions. Although saturation effects and nonlinear interactions are acknowledged in empirical systems, we defer their inclusion to companion studies.

The model tracks the evolution of three interdependent state variables:
\begin{itemize}
    \item $S(t)$: Soil organic matter content (tC/ha), representing long-term soil fertility;
    \item $R(t)$: Residue biomass available at time $t$ (t/ha);
    \item $E(t)$: Cumulative bioenergy output up to time $t$ (MJ/ha).
\end{itemize}
We use $u(t)\in[0,1]$ exclusively as the \emph{allocation} control (fraction of residues diverted to energy). Where $u(t)$ denotes the proportion of biomass directed to energy production, and the remainder $1 - u(t)$ is returned to the soil.

The system is governed by the following set of coupled differential equations:
\begin{align}
    \dot{S}(t) &= \alpha (1 - u(t)) R(t) - \delta_S S(t) + \theta E(t), \label{eq:soil} \\
    \dot{R}(t) &= \eta(S) - \gamma R(t), \label{eq:residue} \\
    \dot{E}(t) &= \beta u(t) R(t) - \delta_E E(t), \label{eq:energy}
\end{align}
with initial conditions:
\begin{equation}
    S(0) = S_0 > 0, \quad R(0) = R_0 > 0, \quad E(0) = 0. \label{eq:initial}
\end{equation}

\textbf{Each equation models a distinct dynamic process.} Equation~\eqref{eq:soil} describes the evolution of soil fertility as a balance between enrichment from biomass return, natural decay, and energy-driven reinforcement. Equation~\eqref{eq:residue} links residue biomass availability to soil quality and includes losses due to decomposition or diversion. We assume residue productivity is linear: $\eta(S)=\rho S$ with $\rho>0$.
 Equation~\eqref{eq:energy} captures energy accumulation through conversion, with degradation of usable output over time.

The model depends on the following parameters:
\begin{itemize}
    \item $\alpha$ [tC/t]: Efficiency of soil enrichment via residue return;
    \item $\delta_S$ [1/year]: Natural decay rate of organic matter in the soil;
    \item $\rho$ [t/(tC$\cdot$year)]: Residue productivity per unit soil fertility;
    \item $\gamma$ [1/year]: Loss rate of residue biomass;
    \item $\beta$ [MJ/t]: Energy yield per ton of residue;
    \item $\delta_E$ [1/year]: Decay rate of accumulated energy;
    \item $\theta$ [tC/(MJ$\cdot$year)]: Efficiency of energy reinvestment into soil.
\end{itemize}

\paragraph*{Modeling rationale (linearity).}
We adopt a linear–affine structure as a calibrated local model around the operating corridor of interest. On the soil side, first-order input–decay kinetics are widely used for SOC turnover and accurately fit long-term bare-fallow datasets across diverse edaphic and climatic conditions \cite{menichetti2019}. On the agronomic side, crop/soil–productivity responses are well described by concave, saturating laws (e.g., Mitscherlich “diminishing returns”), for which a linear law is the first-order approximation near the calibrated range \cite{dhanoa2022}. Within this corridor, the linear–affine specification yields parsimony, identifiability, and transparent optimality conditions.

\begin{remark}
The inclusion of $E(t)$ as a dynamic state variable plays a pivotal role in shaping system behavior. By linking past energy production to current soil dynamics, the model introduces a form of delayed reinforcement that affects both the timing and structure of optimal strategies. This feature supports simulation of sustainability-oriented policies such as reinvestment subsidies or incentive-based conservation programs.
\end{remark}

\begin{lemma}[Positivity] 
\label{lemma:positivity}
Let $u(\cdot)\in[0,1]$, $\alpha,\beta,\gamma,\delta_S,\delta_E>0$, and $\theta\ge 0$. If $S(0),R(0),E(0)\ge 0$, then the solution of \eqref{eq:soil}–\eqref{eq:energy} satisfies $S(t),R(t),E(t)\ge 0$ for all $t\in[0,T]$.
\end{lemma}
\begin{proof}
From \eqref{eq:energy}, $\dot E=\beta uR-\delta_E E\ge -\delta_E E$ with $E(0)\ge 0$. By comparison/Gronwall, $E(t)\ge 0$ for all $t$ (standard Grönwall estimate). 
From \eqref{eq:residue}, $\dot R=\rho S-\gamma R\ge -\gamma R$ and $R(0)\ge 0$ imply $R(t)\ge 0$. 
From \eqref{eq:soil}, $\dot S=\alpha(1-u)R-\delta_S S+\theta E\ge -\delta_S S$ and $S(0)\ge 0$ yield $S(t)\ge 0$. 
Hence the nonnegative orthant is forward positivity for the system.\qedhere
\end{proof}

\begin{remark}
    A companion manuscript develops a multi-objective variant with explicit fertility constraints and mild saturation in the soil–productivity channel, showing that the qualitative bang–interior–bang policy structure persists under those extensions \cite{dekkakiUnderReview}. We therefore retain the linear–affine baseline here and refer readers to that analysis for the saturated case.
\end{remark}

\section{Statement and Existence of the Optimal Control Problem}
\label{sec:ocp}

The goal is to determine a time-dependent control policy $u^*(\cdot)$ that maximizes the total discounted net benefit over a finite horizon $[0, T]$. The benefit reflects a trade-off between energy revenue and soil fertility value, penalized by operational costs related to biomass logistics and energy conversion.

The objective functional is given by:
\begin{equation}
    J(u) = \int_0^T e^{-\delta t} \left[ p_E \dot{E}(t) + p_S S(t) - C(u(t)) \right] \, dt,
    \label{eq:objective}
\end{equation}
where:
\begin{itemize}
    \item $\dot{E}(t)$ [GJ$\cdot$ha$^{-1} \cdot$yr$^{-1}$]: instantaneous energy output rate;
    \item $S(t)$ [tC$\cdot$ha$^{-1}$]: soil organic matter, representing long-term soil fertility;
    \item $p_E$ [\$/MJ]: unit price of usable bioenergy;
    \item $p_S$ [\$/tC]: economic valuation of soil fertility;
    \item $C(u) = c_1 u + c_2 u^2$ [\$/ha$\cdot$yr]: operational cost of diverting residues to energy, with:
    \begin{itemize}
        \item $c_1$ [\$/ha$\cdot$yr]: linear cost coefficient (e.g., residue collection, transport);
        \item $c_2$ [\$/ha$\cdot$yr]: nonlinear cost coefficient (e.g., processing complexity, diminishing efficiency);
    \end{itemize}
    \item $\delta$ [yr$^{-1}$]: economic or social discount rate.
\end{itemize}

The control function $u(t)$ is measurable and constrained to:
\begin{equation}
    \mathcal{U} := \left\{ u(\cdot) \in L^\infty([0, T]; [0, 1]) \right\},
\end{equation}
with $u(t) = 0$ indicating full residue return to soil and $u(t) = 1$ indicating complete diversion to energy.

\begin{remark}
The quadratic cost function $C(u) = c_1 u + c_2 u^2$ reflects both basic operational expenses and additional costs that grow with higher diversion rates. The convex shape accounts for practical challenges such as reduced processing efficiency, limited residue availability, or increased labor demand when more biomass is sent to energy. This makes the cost structure realistic and grounded in actual logistical trade-offs not just a mathematical simplification.
\end{remark}

The optimal control problem becomes:
\begin{equation*}
    \max_{u(\cdot) \in \mathcal{U}} \quad J(u),
\end{equation*}
subject to the dynamic system \eqref{eq:soil}--\eqref{eq:energy} and initial conditions \eqref{eq:initial}.

Before proceeding with the application of the Pontryagin Maximum Principle, it is essential to establish that the optimal control problem stated above is well-posed. In particular, we must ensure that the admissible control set is nonempty and that the system of differential equations admits bounded solutions for any admissible control. Furthermore, we seek to guarantee the existence of at least one optimal control that maximizes the objective functional over the given time horizon. These foundational properties are formalized in the following theorem.

\begin{theorem}[Existence of Optimal Control and Boundedness of State Trajectories]
\label{:existence}
Fix a final time $T > 0$. Suppose:
\begin{enumerate}
    \item The production law is fixed to the working model $\eta(S)=\rho S$ with $\rho>0$. 
    \item The cost function $C:[0,1]\to\mathbb{R}_+$ is continuous. 
    \item The admissible control set is
    \[
    \mathcal{U} := \left\{ u \in L^\infty([0,T]; [0,1]) \right\};
    \]
    \item The parameters $\alpha, \delta_S, \delta_E, \gamma, \beta, p_E, p_S, \delta$ are strictly positive and $\theta \geq 0$;
    \item The initial state satisfies $S(0) = S_0 > 0$, $R(0) = R_0 > 0$, and $E(0) = 0$.
\end{enumerate}
Then the optimal control problem:
\[
\max_{u \in \mathcal{U}} \quad J(u) := \int_0^T e^{-\delta t} \left[ p_E \dot{E}(t) + p_S S(t) - C(u(t)) \right] dt,
\]
subject to the dynamic system:
\begin{align*}
    \dot{S}(t) &= \alpha (1 - u(t)) R(t) - \delta_S S(t) + \theta E(t), \\
    \dot{R}(t) &= \rho S(t) - \gamma R(t), \\
    \dot{E}(t) &= \beta u(t) R(t) - \delta_E E(t),
\end{align*}
admits at least one optimal control $u^* \in \mathcal{U}$. Moreover, the associated state trajectories $(S^*, R^*, E^*)$ remain uniformly bounded on $[0, T]$.
\end{theorem}

\begin{proof}
\textit{Step 1: Boundedness of state trajectories.}
By the positivity lemma \ref{lemma:positivity}, $S,R,E\ge0$ on $[0,T]$. Hence
\[
\dot S \le \alpha R - \delta_S S + \theta E,\quad
\dot R = \rho S - \gamma R,\quad
\dot E \le \beta R - \delta_E E.
\]
Summing yields
\[
\frac{d}{dt}(S+R+E) \le (\rho-\delta_S)S + (\alpha+\beta-\gamma)R + (\theta-\delta_E)E
 \le C_0\,(S+R+E),
\]
for $C_0 := \max\{\rho-\delta_S,\,\alpha+\beta-\gamma,\,\theta-\delta_E\}+\delta_S+\gamma+\delta_E$. By Grönwall,
\[
S(t)+R(t)+E(t)\le (S_0+R_0+E_0)\,e^{C_0 t}\quad\text{for }t\in[0,T],
\]
so $(S,R,E)$ are uniformly bounded on $[0,T]$.

\textit{Step 2: Existence of an optimal control.}
The right-hand side is Carathéodory and control-affine in $u$; the control set $[0,1]$ is compact, so the velocity set $f(t,x,[0,1])$ is convex. Moreover, substituting $\dot E=\beta uR-\delta_E E$ gives a continuous running payoff
\[
\ell(t,x,u)=e^{-\delta t}\big(p_E(\beta uR-\delta_E E)+p_S S - C(u)\big),
\]
continuous in $(t,x,u)$ on the compact control set. By the Filippov–Cesari existence theorem for Bolza problems with compact controls and convex velocity sets (see, e.g., Clarke,\cite{clarke13}; Vinter, \cite{vinter2000}), there exists an optimal $u^*\in\mathcal U$ with associated bounded state $(S^*,R^*,E^*)$.
\end{proof}

\section{Pontryagin Maximum Principle Analysis}
\label{sec:PMP}

To derive necessary conditions for optimality in the crop residue allocation problem, we apply the Pontryagin Maximum Principle (PMP) in the \emph{current-value} formulation. This approach simplifies the structure of the Hamiltonian by eliminating the discounting exponential $e^{-\delta t}$ from the integrand in the objective functional \eqref{eq:objective}, making both analytical and numerical treatments more tractable. The structure follows the methodology presented in \cite{kamien2012dynamic}, Sect. 8 of Part II.

Let $\lambda = (\lambda_S, \lambda_R, \lambda_E)$ be the vector of costate (adjoint) variables associated with the state variables $(S, R, E)$. We use the term \emph{current-value pseudo-costate} for $e^{\delta t}\lambda_i(t)$, $i\in\{S,R,E\}$, where $\lambda_i$ are the standard costates in current-value form.
 The current-value Hamiltonian for our system is defined as:
\begin{equation}
    H = p_E \beta u R + p_S S - C(u) + \lambda_S[\alpha(1 - u) R - \delta_S S+ \theta E] + \lambda_R[\eta(S) - \gamma R] + \lambda_E[\beta u R - \delta_E E].
    \label{eq:hamiltonian}
\end{equation}

We reorganize the Hamiltonian as follows:
\begin{equation}
    H = h(S, R, \lambda) + \tilde{h}(\lambda) u R - C(u),
\end{equation}
where:
\begin{itemize}
    \item $h(S, R, \lambda) := p_S S + \lambda_S \alpha R - \lambda_S \delta_S S + \lambda_S \theta E + \lambda_R [\eta(S) - \gamma R]- \lambda_E \delta_E E$,
    \item $\tilde{h}(\lambda) := p_E \beta - \lambda_S \alpha + \lambda_E \beta$.
\end{itemize}

The adjoint dynamics are given by:
\begin{equation}
    \left\{
    \begin{array}{lll}
        \dot{\lambda}_S &= \delta \lambda_S - \dfrac{\partial H}{\partial S} = \delta \lambda_S - \left[p_S - \lambda_S \delta_S + \lambda_R \eta'(S) \right], \\[1.5mm]
        \dot{\lambda}_R &= \delta \lambda_R - \dfrac{\partial H}{\partial R} = \delta \lambda_R - \left[ (p_E + \lambda_E) \beta u - \lambda_S \alpha(1 - u) - \lambda_R \gamma \right], \\[1.5mm]
        \dot{\lambda}_E &= \delta \lambda_E - \left[ \theta \lambda_S - \lambda_E \delta_E \right].
    \end{array}
    \right.
    \label{eq:costates}
\end{equation}

Since the final state is free, the transversality conditions are:
\begin{equation}
    \lambda_S(T) = \lambda_R(T) = \lambda_E(T) = 0.
    \label{eq:transversality}
\end{equation}

\textbf{Switching function:}
\begin{equation}
  \phi(t)\,=\,\frac{\partial H}{\partial u}(t)
  \;=\;\big(p_E\beta-\alpha\,\lambda_S(t)+\beta\,\lambda_E(t)\big)\,R(t)\;-\;\big(c_1+2c_2\,u(t)\big).
  \label{eq:switching}
\end{equation}
Bang at bounds when $\phi(t)\lessgtr 0$; interior $u(t)$ solves $\phi(t)=0$.

The PMP maximization condition implies that for almost every $t \in [0,T]$,
\begin{equation}
    u^*(t) \in \arg\max_{u \in [0,1]} H(t).
\end{equation}

To determine the explicit form of the optimal control $u^*(t)$, we apply the Karush-Kuhn-Tucker (KKT) conditions directly, using the fact that the control cost is quadratic and the control is bounded.

\begin{proposition}[Structure of the Optimal Control]
\label{prop:optcontrol}
Let $u^*(t)$ denote the optimal control for the problem defined in Section~\ref{sec:model}, and suppose $C(u) = c_1 u + c_2 u^2$ with $c_1, c_2 > 0$. Then for almost every $t \in [0,T]$, $u^*(t)$ is given by:
\begin{equation}
    u^*(t) =
    \begin{cases}
        0, & \text{if } \tilde{h}(t) R(t) \leq c_1, \\[1mm]
        \min\left\{ \dfrac{\tilde{h}(t) R(t) - c_1}{2 c_2}, 1 \right\}, & \text{if } \tilde{h}(t) R(t) > c_1,
    \end{cases}
    \label{eq:optcontrol}
\end{equation}
where $\tilde{h}(t) := p_E \beta - \lambda_S(t) \alpha + \lambda_E(t) \beta$.
\end{proposition}

\begin{proof}
Let $L(u, W_1, W_2)$ denote the Lagrangian associated with the constrained maximization of the Hamiltonian:
\begin{equation*}
    L = H + W_1 (0 - u) + W_2 (u - 1), \quad \text{with } W_1, W_2 \geq 0.
\end{equation*}
The KKT conditions state that the optimal control $u^*(t)$ satisfies:
\begin{align*}
    &\frac{\partial L}{\partial u} = \frac{\partial H}{\partial u} - W_1 + W_2 = 0, \\[1mm]
    &W_1 u^*(t) = 0, \quad W_2 (1 - u^*(t)) = 0.
\end{align*}
Since $H$ is differentiable and $C(u) = c_1 u + c_2 u^2$, we have:
\begin{equation*}
    \frac{\partial H}{\partial u} = \tilde{h}(t) R(t) - C'(u(t)) = \tilde{h}(t) R(t) - (c_1 + 2c_2 u(t)).
\end{equation*}
Setting this derivative equal to zero for interior solutions yields:
\begin{equation*}
    u(t) = \frac{\tilde{h}(t) R(t) - c_1}{2 c_2}.
\end{equation*}
We then enforce the box constraints $u \in [0,1]$ using complementary slackness:
\begin{equation*}
    u^*(t) = \min \left\{ \max \left\{ 0, \frac{\tilde{h}(t) R(t) - c_1}{2 c_2} \right\}, 1 \right\}.
\end{equation*}
This expression leads to the piecewise structure in \eqref{eq:optcontrol}.
\end{proof}

The value of $\tilde{h}(t)$ determines the marginal benefit of allocating crop residues to bioenergy. If this marginal gain is too low, all residue is returned to soil; if it is high, full diversion to energy occurs. Intermediate values yield interior control values.

\paragraph{Link to application.}
In the crop–energy setting, the current-value costates $\lambda_S,\lambda_R,\lambda_E$ are the shadow values of soil fertility, residue stock, and accumulated energy. The switching function $\phi(t)=\partial H/\partial u$ compares the marginal benefit of sending one more unit of residue to energy with its opportunity cost for soil: $\phi(t){>}0$ pushes $u$ to the upper bound (more energy), $\phi(t){<}0$ pushes $u$ to the lower bound (more soil), and $\phi(t)\approx0$ yields an interior allocation where marginal values balance.
This decision rule is the operational link from PMP theory to residue-management policy (how much to divert vs.\ return at each time). Over medium/long horizons, the observed mid-horizon interior phase acts as a quasi-turnpike: the system stays near a steady allocation for most of the horizon, moving only near the endpoints.
Relative to our previous waste-to-energy OCPs, which establish the PMP framework and direct-method verification on waste recovery and investment,\cite{cherkaoui2022bio, cherkaoui2023acc} the present model adds a reinforcement channel from accumulated energy back to soil (reinvestment), creating the explicit energy–soil feedback that drives the interior phase and the comparative statics in Section~\ref{sec:numerics}.

\section{Numerical Simulations}
\label{sec:numerics}

To approximate the solution of the optimal control problem formulated in Section~\ref{sec:PMP}, we implement a direct optimization approach using a discretization-based method. This approach transforms the continuous-time optimal control problem into a nonlinear programming (NLP) problem, which can then be solved efficiently using modern solvers.

All simulations are performed in the {\tt Julia} programming language, leveraging the {\tt JuMP} modeling interface and the interior-point solver {\tt Ipopt}. Following the methodology outlined in \cite{jbcaillau22}, we discretize the state and control trajectories using a {\tt Crank-Nicolson} scheme, which provides a good balance between accuracy and stability.

The numerical resolution proceeds by dividing the time interval $[0,T]$ into $N$ uniform steps. The control variable $u(t)$ is approximated as piecewise constant over each time step, and the state variables $(S(t), R(t), E(t))$ are approximated using implicit midpoint integration. This allows the resulting NLP to preserve the smooth structure of the continuous system and support a high convergence precision.

The discretization settings and solver tolerances used in the simulations are summarized in Table~\ref{Tab.disc}. The values for the model parameters and the economic coefficients used in the objective function are provided in Table~\ref{Tab.param}. These values are selected to reflect plausible trade-offs between soil fertility and bioenergy production, consistent with empirical studies in residue management and bioenergy systems.

\begin{table}[h!]
	\centering
	\caption{Discretization and solver configuration settings \label{Tab.disc}}
	\begin{tabular}{l c}
		\hline
		\textbf{Discretization method} & {\tt Crank-Nicolson} \\
		\hline
		Number of time steps ($N$) & $2000$ \\
		\hline
		NLP tolerance (absolute) & $10^{-10}$ \\
		\hline
		Control discretization & Piecewise constant \\
		\hline
	\end{tabular}
\end{table}

\begin{table}[h!]
\centering
\caption{Model Parameters: Values and Interpretations \label{Tab.param}}
\begin{tabular}{llp{9cm}}
\toprule
\textbf{Parameter} & \textbf{Value} & \textbf{Interpretation}\\
\midrule
$\gamma$ (Residue decay rate) & 0.20 & Share of residues lost naturally over time; \textit{estimated from seasonal soil decay data}. \\
$\delta_S$ (Soil degradation rate) & 0.05 & Annual decline in soil fertility without residue input; \textit{assumed low to reflect slow degradation}. \\
$\alpha$ (Residue return efficiency) & 0.25 & Proportion of returned biomass that boosts soil fertility; \textit{calibrated to represent moderate composting efficiency}. \\
$\beta$ (Energy conversion efficiency) & 0.35 & Energy yield from converting residues; \textit{estimated from anaerobic digestion efficiency ranges}. \\
$\rho$ (Soil productivity factor) & 0.50 & Contribution of soil fertility to crop residue production; \textit{assumed for moderate responsiveness}. \\
$\delta$ (Discount rate) & 0.02 & Reflects economic preference for present vs future benefits; \textit{standard assumed annual rate}. \\
$p_E$ (Unit energy price) & 1.00 & Normalized unit price of bioenergy; \textit{baseline for comparative analysis}. \\
$p_S$ (Soil fertility value) & 0.80 & Economic value of soil health per unit; \textit{assumed to weigh against short-term energy gain}. \\
$c_1$, $c_2$ (Cost coefficients) & 0.80, 1.00 & Marginal and quadratic costs of diverting biomass to energy; \textit{calibrated to yield interior control}. \\
$T$ (Planning horizon) & 25 & Total time in years for optimization; \textit{assumed to capture both transient and long-term effects}. \\
$\theta$ (Reinvestment efficiency) & 0.20 & Proportion of energy returns reinvested in soil fertility; \textit{estimated to model circular bioeconomy}. \\
$\delta_E$ (Energy degradation rate) & 0.03 & Loss rate of accumulated usable energy; \textit{assumed to ensure bounded energy stock}. \\
\bottomrule
\end{tabular}
\end{table}

The initial conditions are assumed as follows:
\[
S(0) = 1.0, \quad R(0) = 0.5, \quad E(0) = 0.
\]

\paragraph{Parameter context.}
The normalized energy price $P_E{=}1.0$ is mapped to representative EU anchors: recent feed-in reference prices for \emph{biomass}  \texteuro 153 – \texteuro 176/MWh and \emph{biogas}  \texteuro 192– \texteuro 219/MWh (Greece schedules), and the EU emergency cap on inframarginal market revenues at \texteuro180/MWh.\cite{cleanEUislands,eurcap2022}
A policy-facing proxy for soil valuation $P_S$ is the carbon price: EU ETS reports place 2023–24 averages roughly   \texteuro 65-  \texteuro 85/$tCO_2$ ($\approx$ \texteuro 240 – \texteuro 310 per $tC$ using $1 tC = 3.667$ $tCO_2$).(see, \cite{esma2024,reutersETS2024})

Finally, $\theta=0.2$ is a low–moderate reinforcement consistent with observed soil organic carbon (SOC) responses: meta-analyses report SOC increases from residue return around $\sim$11\% and large, persistent SOC gains under biochar additions.\cite{wang2020,beillouin2023}

\begin{figure}[h!]
  \centering
  \includegraphics[width=\textwidth]{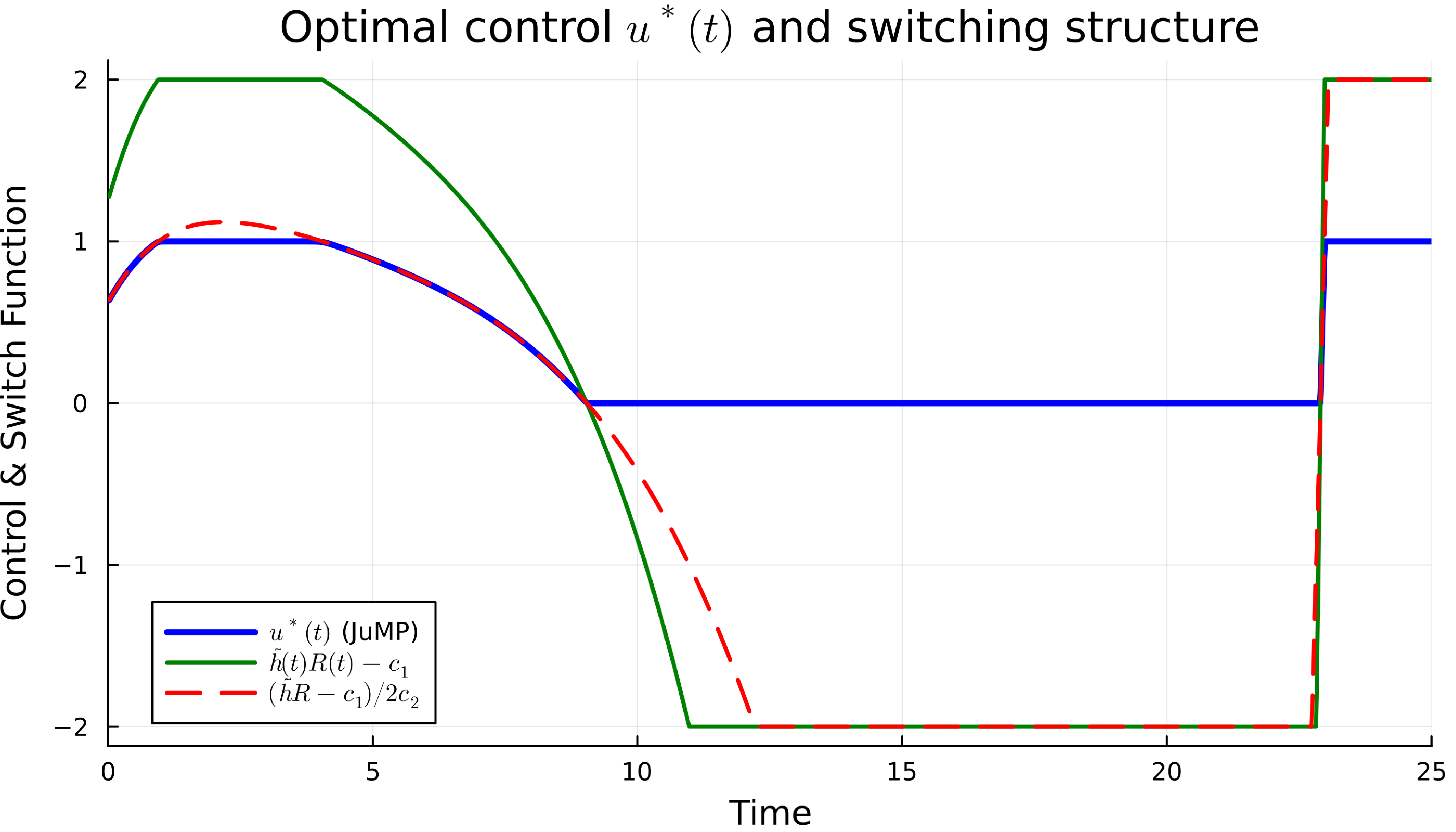}
 \caption{Optimal allocation $u(t)$ (top) and switching function $\phi(t)$ (bottom), cf.\ Eq.~\eqref{eq:switching}. 
Zeros of $\phi$ identify candidate switches; on intervals where $\phi(t)=0$ the interior control solves the KKT condition, 
while $\phi(t){<}0$ (resp.\ $>{}0$) pushes the optimum to the lower (resp.\ upper) bound. Parameters as in Table~\ref{Tab.param}.}
  \label{fig:control_switch}
\end{figure}

\begin{figure}[h!]
  \centering
  \includegraphics[width=\textwidth]{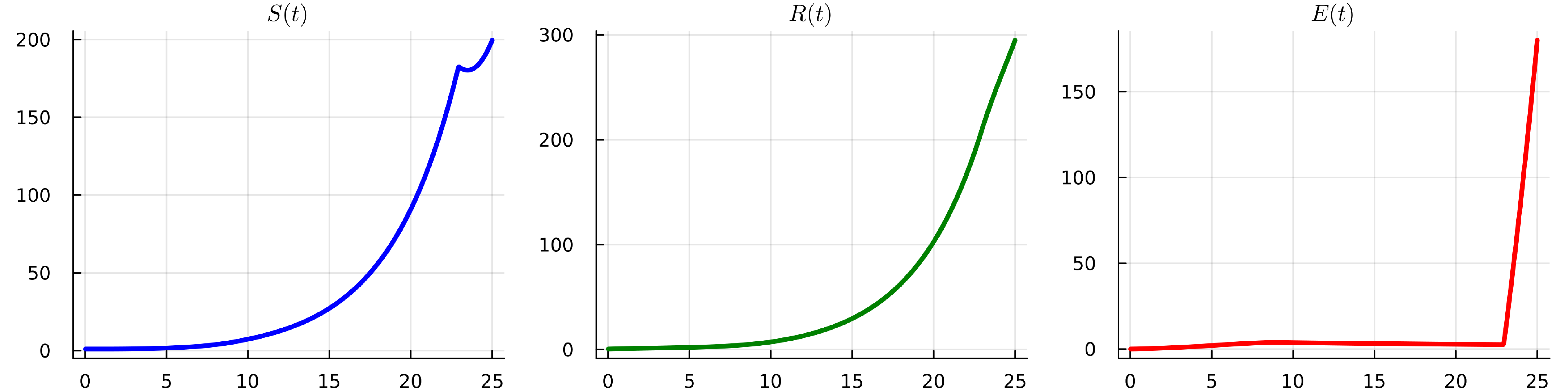}
\caption{States $S(t)$ (soil), $R(t)$ (residue), and $E(t)$ (cumulative energy). 
The mid-horizon plateau coincides with interior operation in $u(t)$; early/late drifts reflect bang phases. 
Parameters as in Table~\ref{Tab.param}.}
  \label{fig:states}
\end{figure}

\begin{figure}[h!]
    \centering
    \includegraphics[width=0.49\textwidth]{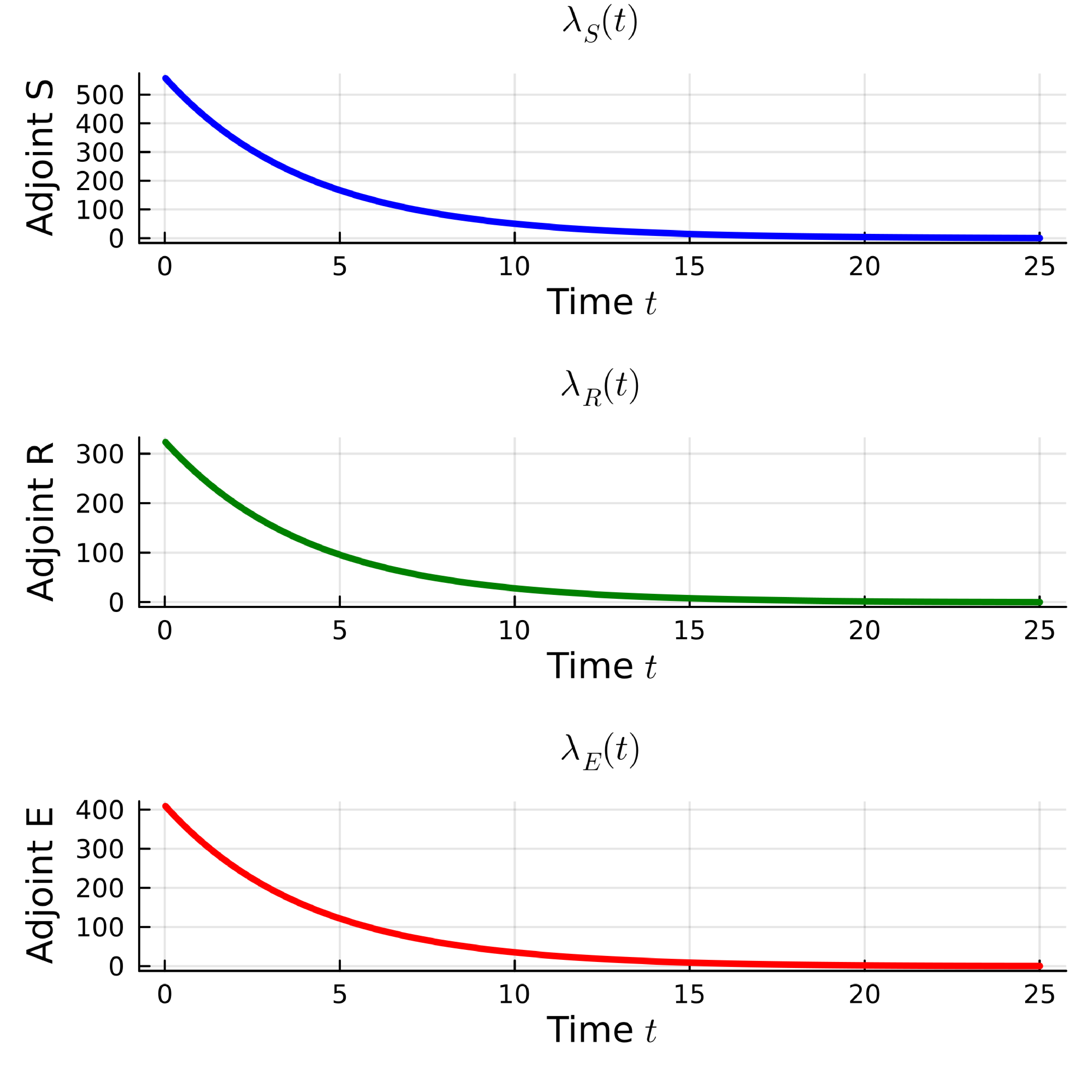}
    \includegraphics[width=0.49\textwidth]{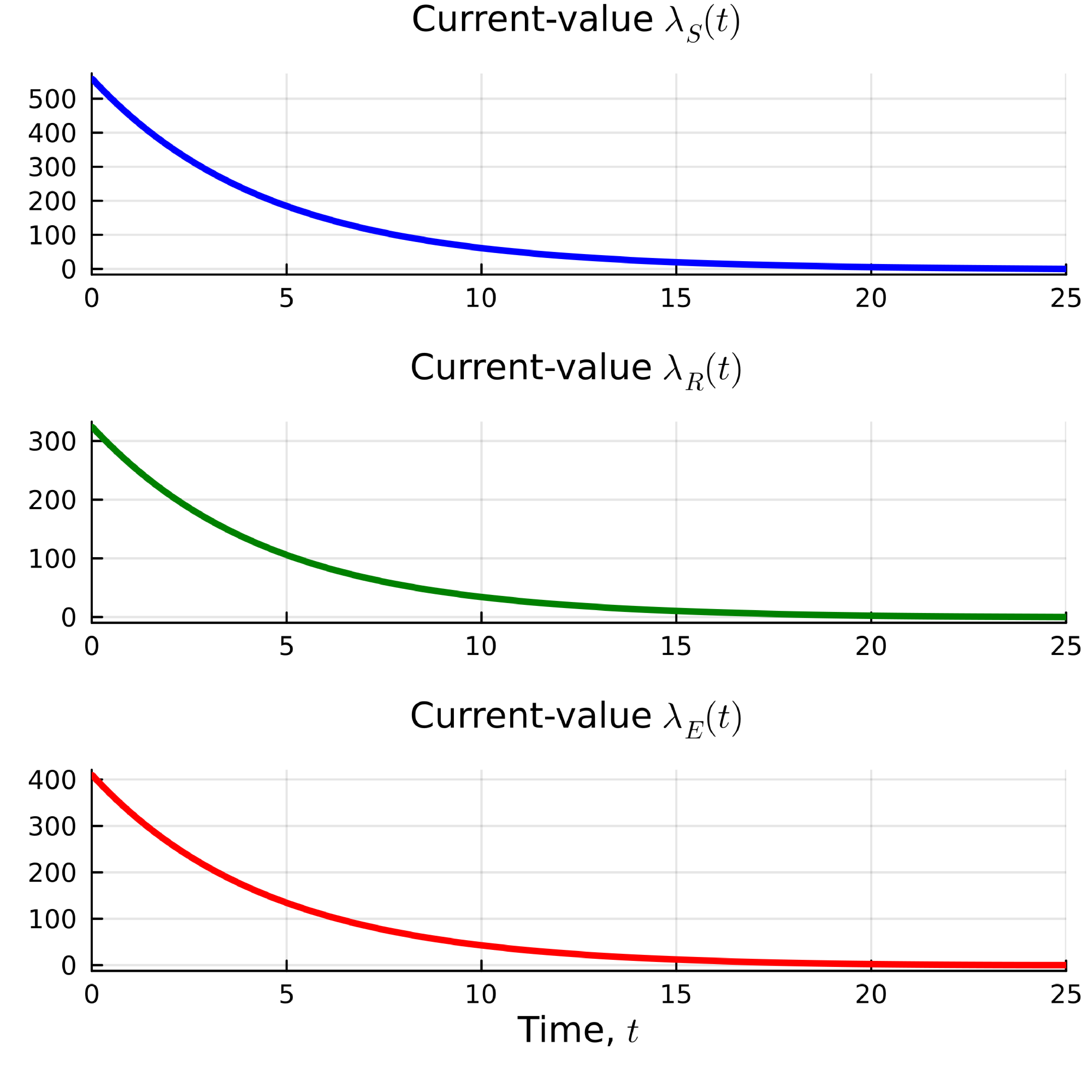}
    \caption{On the left: the adjoint costate trajectories $\lambda_S(t)$, $\lambda_R(t)$, and $\lambda_E(t)$ obtained via {\tt JuMP} for the reinvestment model. On the right: the corresponding current-value pseudo-costates $e^{\delta t} \lambda_S(t)$, $e^{\delta t} \lambda_R(t)$, and $e^{\delta t} \lambda_E(t)$. The decaying profiles reflect the diminishing marginal value of state variables under discounting and confirm consistency with the Pontryagin Maximum Principle.}
    \label{fig:costate_combined}
\end{figure}

The benchmark simulation over a 25-year horizon provides insight into the structure of the optimal policy and the system’s dynamic response to reinvestment and decay mechanisms. The optimal control trajectory and switching structure are shown in Figure~\ref{fig:control_switch}\footnote{To improve visual clarity in Figure~\ref{fig:control_switch}, the switching function $\tilde{h}(t)R(t) - c_1$ and its scaled version were clamped using {\tt clamp01(x) = clamp(x, -2.0, 2.0)}. This operation does not alter the computed control $u^*(t)$, but limits the plotted range of the switching expressions to the interval $[-2, 2]$ for graphical readability.}
, while the corresponding current-value costate variables and adjoint trajectories are reported in Figure~\ref{fig:costate_combined}, respectively. The evolution of the state variables $S(t)$, $R(t)$, and $E(t)$ is depicted in Figure~\ref{fig:states}.

The optimal control $u^*(t)$ (blue curve in Figure~\ref{fig:control_switch}) exhibits a nontrivial structure with three distinct regimes. It starts at its upper bound, allocating all residues to energy production in the early phase (up to $t \approx 2.5$), then transitions to an interior control value around $u \approx 0.65$ during a middle interval where reinvestment builds soil fertility while maintaining a moderate energy supply. Around $t \approx 10$, the control switches sharply to its lower bound $u = 0$, fully diverting residues back to the soil. This reflects the shift in the sign of the switching function $\tilde{h}(t)R(t) - c_1$ (green curve), which falls below zero, indicating that the marginal value of energy production is lower than the combined cost and fertility opportunity. As the horizon approaches its end ($t \approx 23$), the control switches again to $u = 1$, favoring immediate energy gain, consistent with the terminal transversality conditions and time-discounted preferences.

The behavior of the adjoint variables, shown in Figure~\ref{fig:costate_combined}, confirms the structure predicted by the Pontryagin Maximum Principle. Each adjoint $\lambda_i(t)$ decays monotonically toward zero, satisfying the transversality condition $\lambda_i(T) = 0$ for $i \in \{S, R, E\}$. The initial values of $\lambda_S(0)$ and $\lambda_E(0)$ are significantly higher than that of $\lambda_R(0)$, reflecting the high marginal shadow prices of long-term soil fertility and usable energy early in the planning horizon.

State trajectories (Figure~\ref{fig:states}) show rapid exponential growth of $S(t)$ and $R(t)$, consistent with the positive feedback created by reinvestment. The energy stock $E(t)$ remains nearly flat over most of the horizon but rises sharply in the final phase when $u(t) \to 1$ again. This confirms that the system initially prioritizes building capital (fertility and productivity), then exploits the accumulated potential near the horizon.

Overall, the optimal solution highlights a dynamically efficient policy that balances early reinvestment with late-stage energy exploitation. It confirms that even under energy decay and reinvestment efficiency losses, a three-phase control policy can maximize the cumulative benefit from both soil and bioenergy resources.

\subsection{Comparison with the No-Reinvestment Scenario ($\theta = 0$)}

\begin{figure}[h!]
  \centering
  \includegraphics[width=\textwidth]{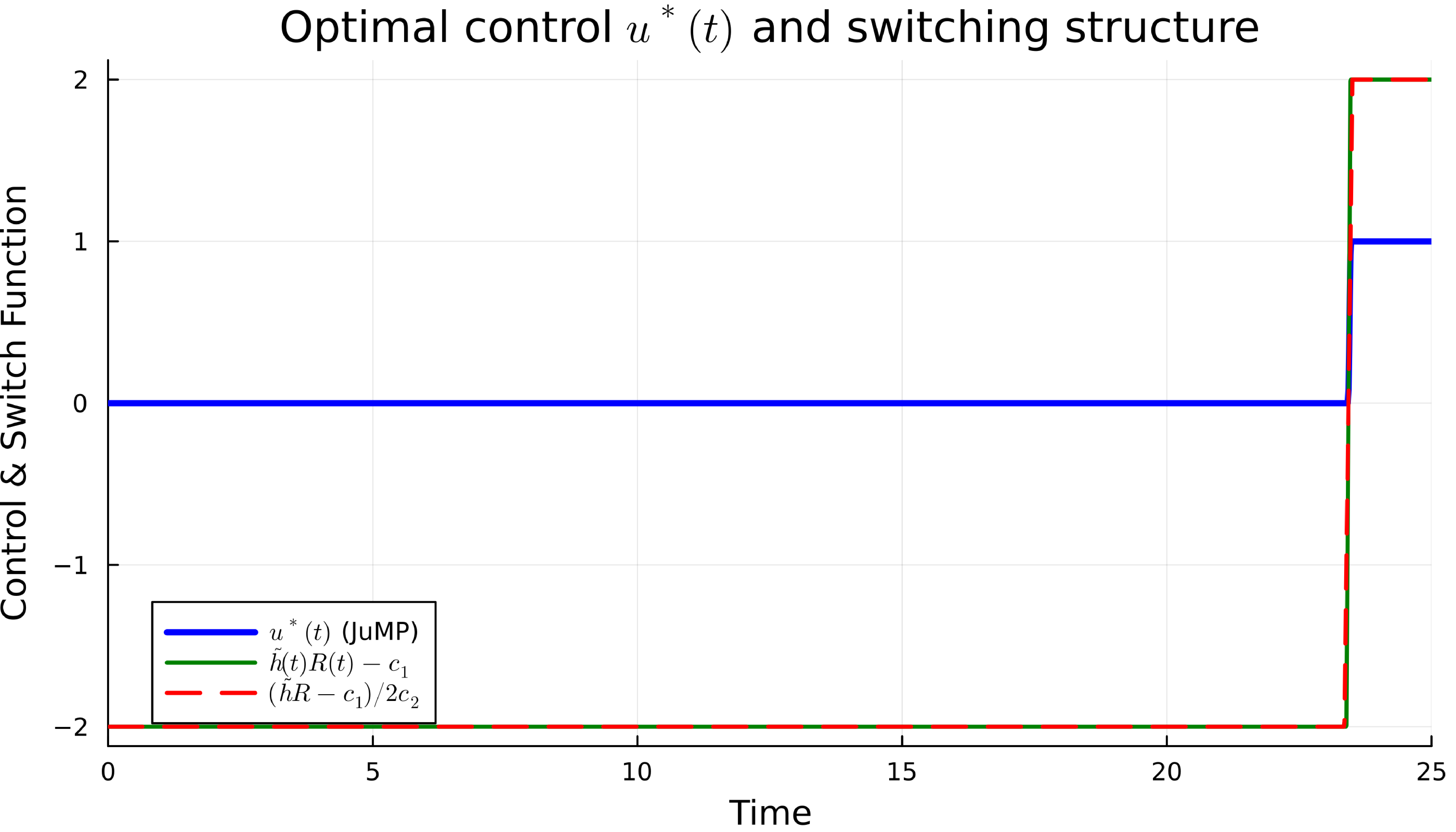}
\caption{No-reinvestment case ($\theta{=}0$): optimal allocation $u(t)$ and switching function $\phi(t){=}\partial H/\partial u$. Zeros of $\phi$ indicate switching times; $\phi{<}0\Rightarrow u{=}0$ (all residues to soil), $\phi{>}0\Rightarrow u{=}1$ (all to energy), and $\phi{\approx}0$ corresponds to interior operation $0{<}u{<}1$. Without $E{\rightarrow}S$ reinforcement, the policy stays at $u{=}0$ for most of the horizon.}
  \label{fig:control_switch_zero}
\end{figure}

\begin{figure}[h!]
  \centering
  \includegraphics[width=\textwidth]{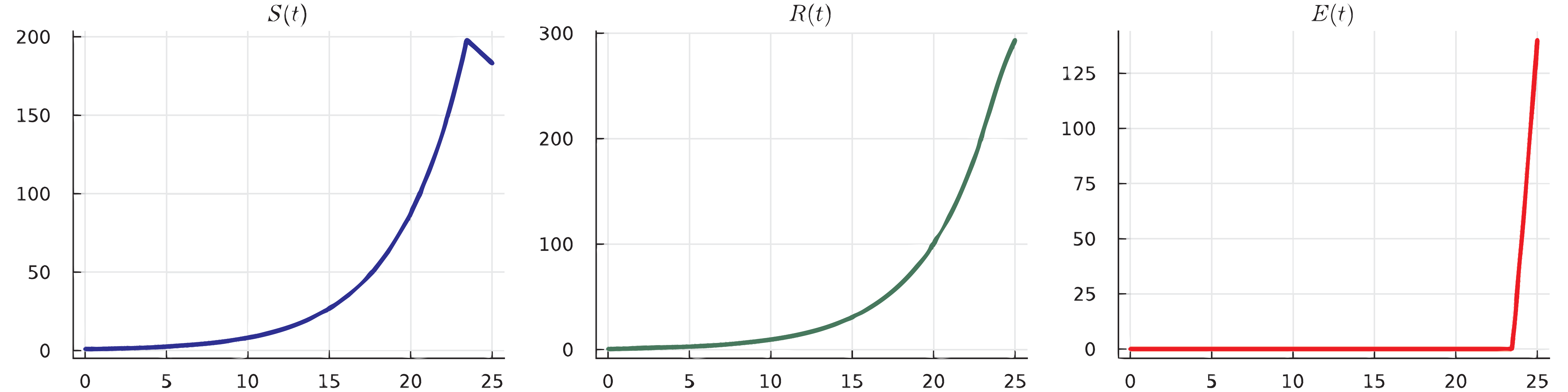}
\caption{No-reinvestment case ($\theta{=}0$): state trajectories $S$ (soil), $R$ (residue), and $E$ (cumulative energy). Reduced growth in $S$ and $R$ reflects sustained $u{=}0$; $E$ rises only near the terminal window where $\phi$ becomes positive.}
  \label{fig:states_zero}
\end{figure}

\begin{figure}[h!]
    \centering
    \includegraphics[width=0.49\textwidth]{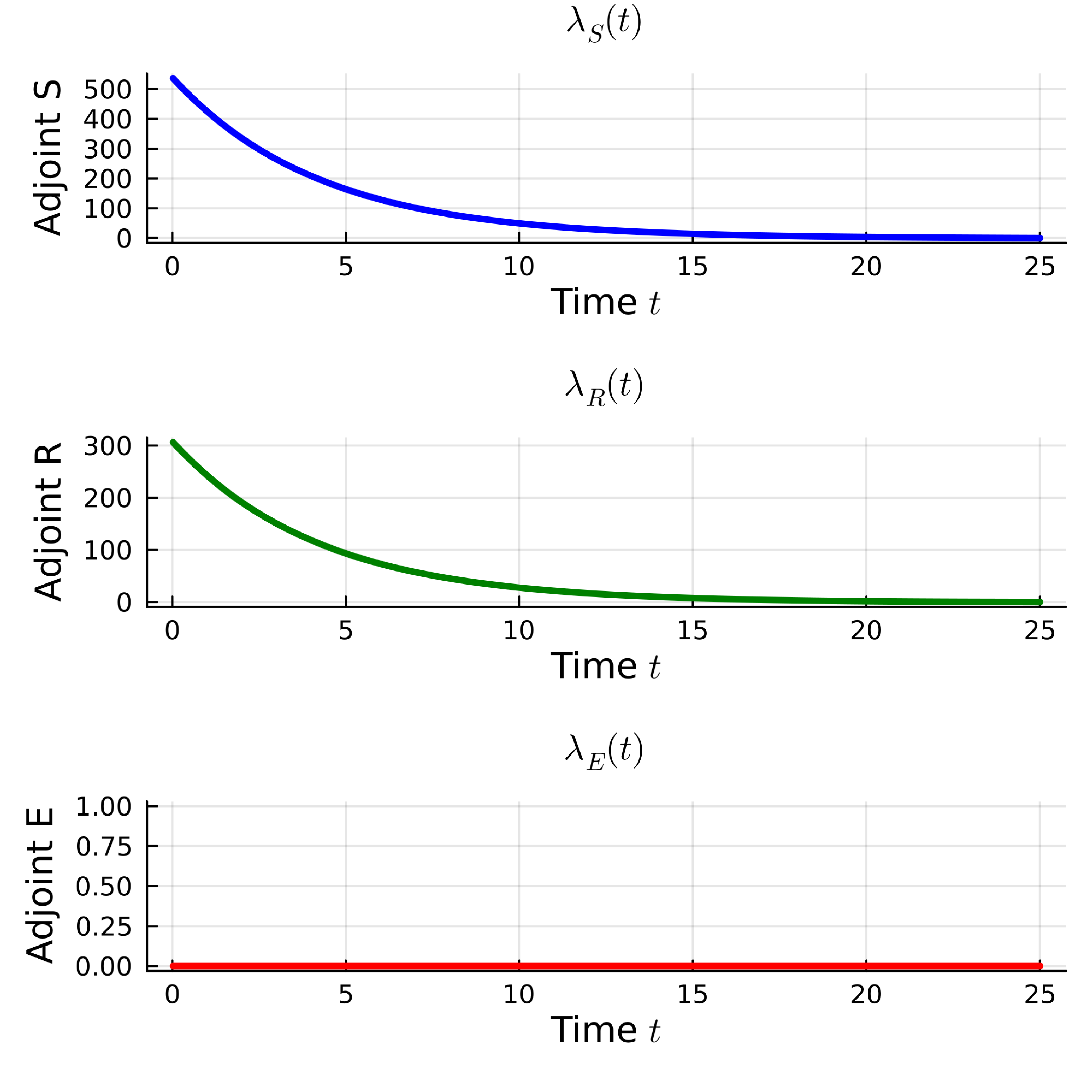}
    \includegraphics[width=0.49\textwidth]{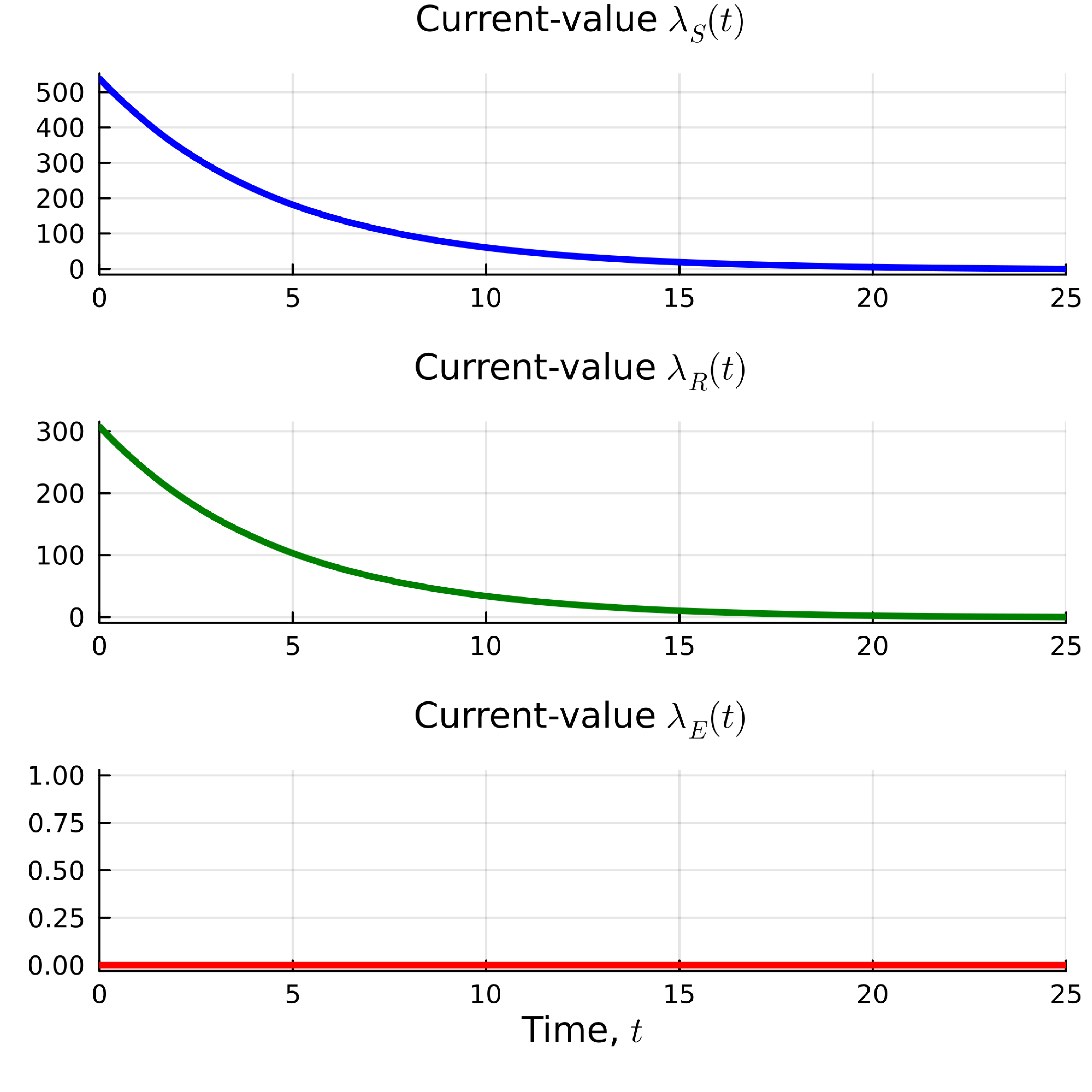}
\caption{No-reinvestment case ($\theta{=}0$): costates $\lambda_S,\lambda_R,\lambda_E$ (left) and current-value pseudo-costates $e^{\delta t}\lambda_i$ (right). Transversality $\lambda_i(T){=}0$ holds; $\lambda_E$ stays near zero, indicating low marginal value of energy without feedback-consistent with $u{=}0$ over most of the horizon.}
    \label{fig:costate_combined_zero}
\end{figure}

To isolate the effect of energy reinvestment on the system's dynamics and the structure of the optimal control policy, we now simulate the model under the same parameters and time horizon, but with the reinvestment coefficient set to zero: $\theta = 0$. In this case, energy accumulation no longer contributes to soil fertility, and the control strategy must rely solely on direct residue return to preserve soil quality. The resulting control trajectory and switching expressions are displayed in Figure~\ref{fig:control_switch_zero}, while the corresponding costate and pseudo-costate trajectories are shown in Figure~\ref{fig:costate_combined_zero}. The resulting state trajectories are plotted in Figure~\ref{fig:states_zero}.

Compared to the reinvestment benchmark, the optimal control $u^*(t)$ (Figure~\ref{fig:control_switch_zero}) displays a markedly different behavior. The control remains pinned at its lower bound $u = 0$ for the entire duration except in the final phase, where a brief surge appears as the switching function crosses the optimality threshold. This indicates that, in the absence of reinvestment, the marginal benefit of energy production is consistently outweighed by the cost and long-term loss in soil fertility.

The costate trajectories in Figures~\ref{fig:costate_combined_zero} confirm this strategic behavior. The current-value costates $\lambda_S(t)$ and $\lambda_R(t)$ dominate the decision structure early on, while $\lambda_E(t)$ remains near zero throughout, indicating that energy has virtually no strategic value when it does not feed back into the soil.

The state dynamics (Figure~\ref{fig:states_zero}) illustrate a significantly weaker trajectory for soil fertility and residue stock compared to the reinvestment scenario. The energy stock $E(t)$ only grows near the end, as the optimizer briefly diverts residue to energy production when the discounting effect renders soil gains less valuable. This contrasts with the benchmark case where reinvestment allowed simultaneous buildup of energy and fertility, supporting a more balanced and efficient control policy.

These results highlight the structural importance of the reinvestment mechanism. Without it, the optimal strategy becomes more conservative and resource-preserving for most of the horizon. Reinvestment not only increases cumulative reward, but also unlocks more aggressive energy policies without compromising sustainability. This reinforces the circular economy perspective introduced in our model.

\subsection{Effect of Planning Horizon on the Optimal Strategy}

\begin{figure}[h!]
    \centering
    \includegraphics[width=0.49\textwidth]{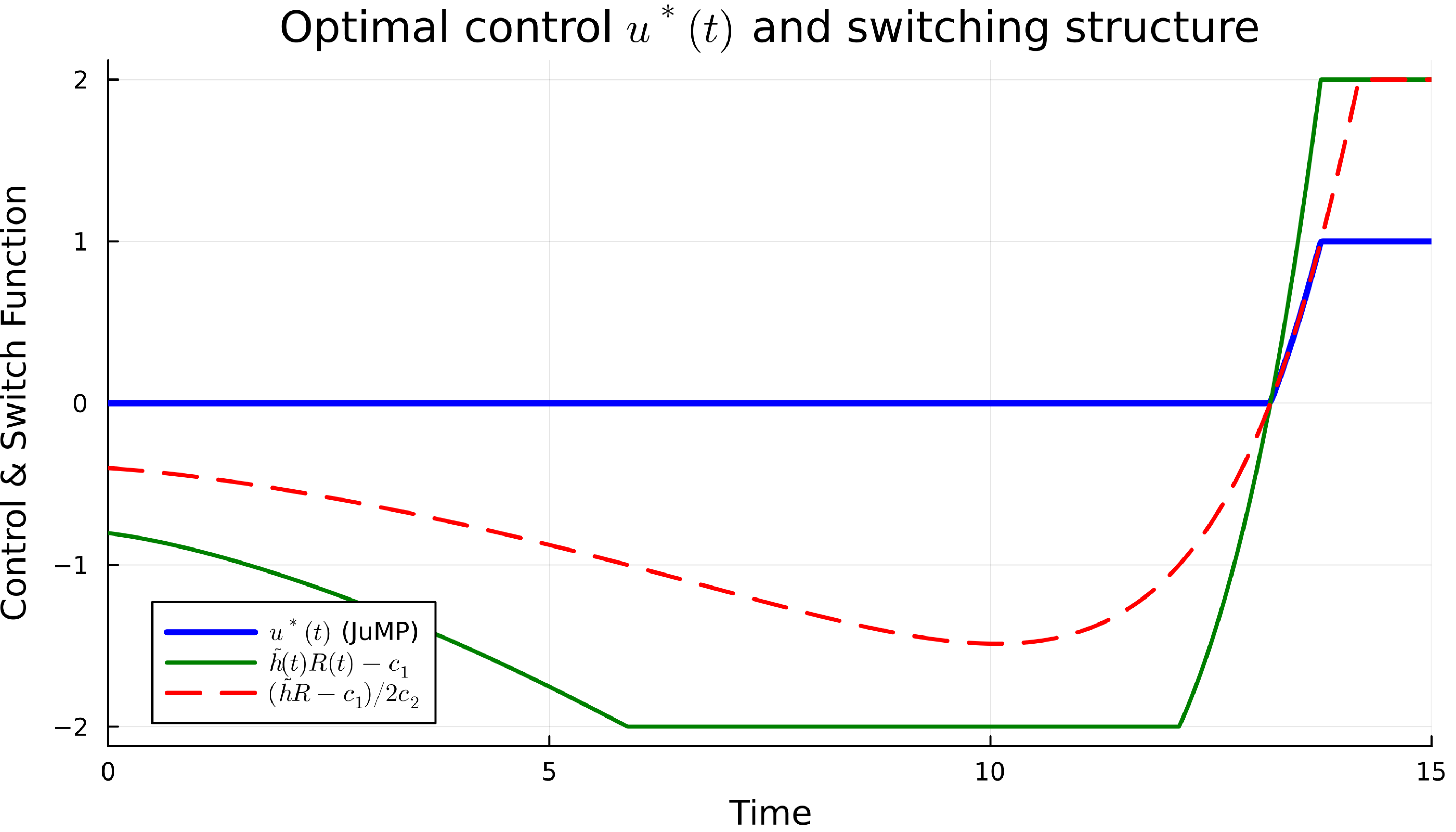}
    \includegraphics[width=0.49\textwidth]{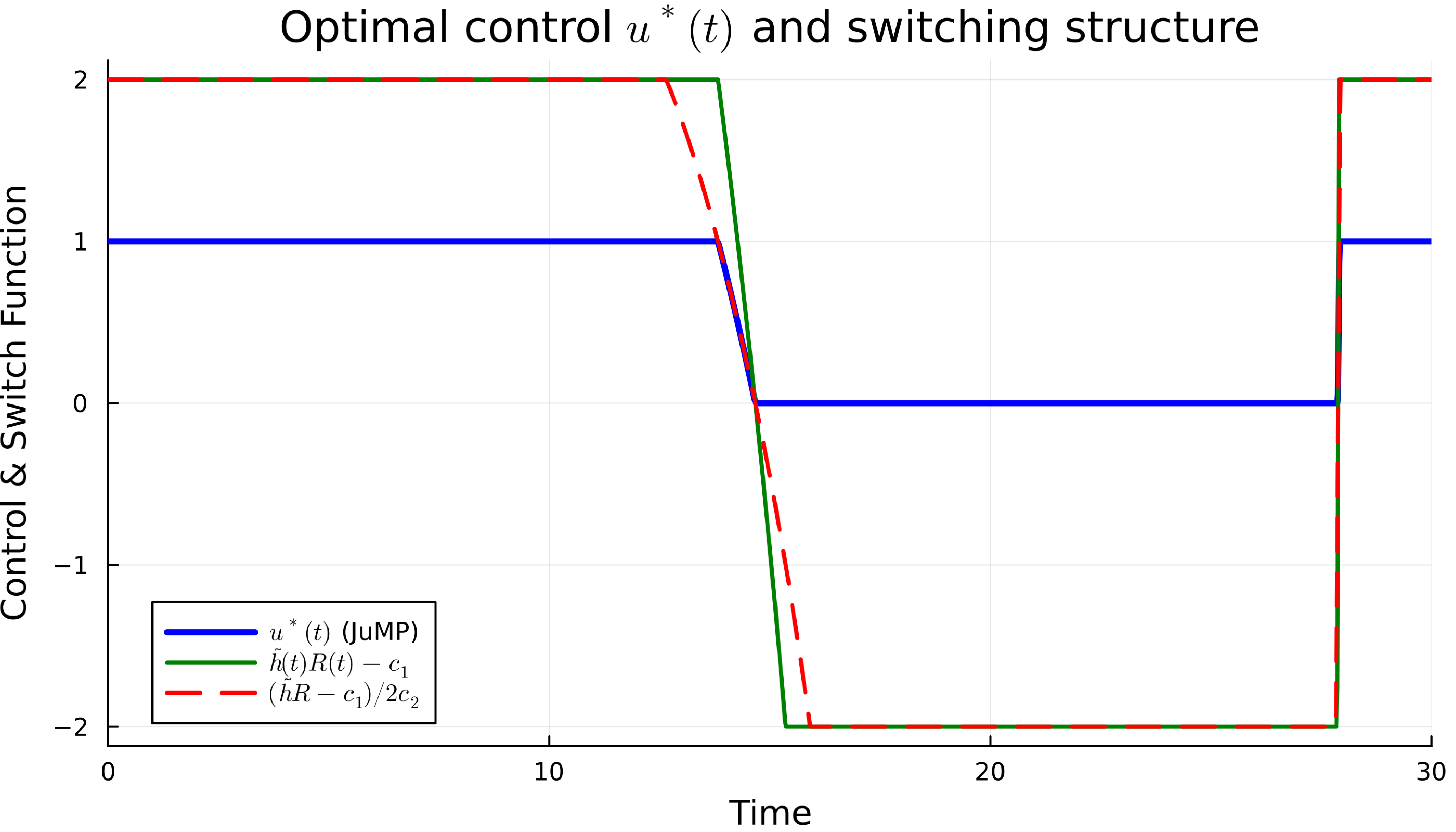}
\caption{Effect of planning horizon on control and switching: short $T{=}15$ (left) vs.\ long $T{=}30$ (right). $u(t)$ is shown with $\phi(t){=}\partial H/\partial u$. Short horizons suppress interior operation (dominant $u{=}0$); long horizons yield a bang–interior–bang structure with a pronounced mid-horizon interior interval.}
    \label{fig:control_horizons}
\end{figure}

\begin{figure}[h!]
    \centering
    \includegraphics[width=0.9\textwidth]{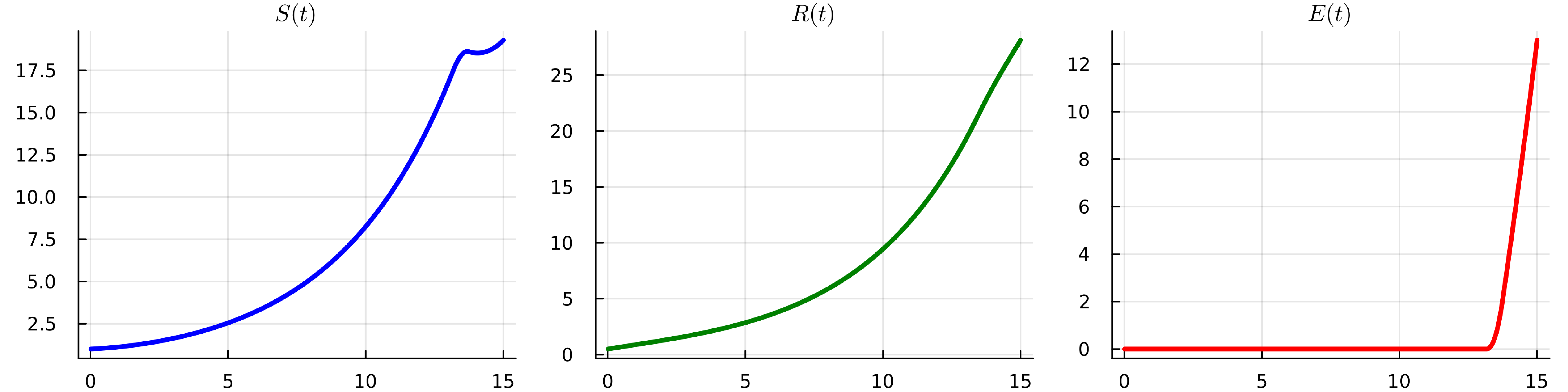}
    \includegraphics[width=0.9\textwidth]{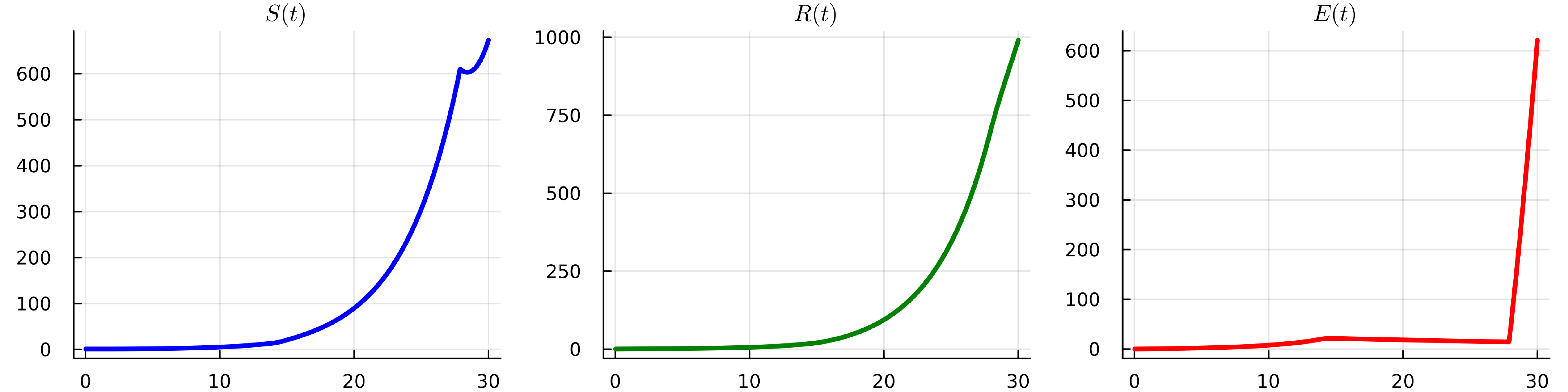}
    \caption{State trajectories under short $T{=}15$ (top) and long $T{=}30$ (bottom) horizons. Longer horizons enable compounding via reinforcement: $S$ and $R$ rise more strongly and $E$ builds earlier, matching the interior phase observed in the control.}
    \label{fig:states_horizons}
\end{figure}

To investigate the impact of the available time horizon on the structure of the optimal control policy and system dynamics, we simulate the model for two scenarios with different final times: a short horizon $T = 15$ and a longer horizon $T = 30$. All parameters, including the reinvestment coefficient $\theta = 0.2$, are kept fixed to isolate the effect of time availability.

Figure~\ref{fig:control_horizons} shows the optimal control trajectories $u^*(t)$ and their switching structure. In the short-horizon case (left panel), the control remains at its lower bound $u = 0$ over most of the time interval, only increasing sharply near the end. This reflects a conservative strategy: there is not enough time to benefit from energy reinvestment, so the optimizer favors residue return to maintain soil fertility. The switching function $\tilde{h}(t)R(t) - c_1$ remains negative for most of the horizon, confirming that energy diversion is not economically justified early on.

In contrast, the long-horizon case (right panel) displays a three-phase structure. The control starts at its upper bound ($u = 1$), then shifts to an interior value, before eventually returning to $u = 0$ and finally switching back to $u = 1$ near the terminal time. This pattern reveals a dynamic reinvestment strategy: early energy production enables mid-horizon soil enhancement, which is later exploited for increased residue and energy output.

The associated state trajectories are shown in Figure~\ref{fig:states_horizons}. In the short-horizon case, fertility $S(t)$ and residues $R(t)$ grow moderately, while energy $E(t)$ rises only in the final phase. In the long-horizon case, all state variables exhibit strong growth. Soil fertility and productivity accumulate over time, supported by reinvestment, and energy builds steadily before accelerating at the end.

These simulations highlight the central role of the planning horizon in shaping the control policy. Short horizons favor conservative allocation, focusing on immediate soil preservation. Long horizons allow more complex strategies involving reinvestment and exploitation phases, illustrating the temporal trade-offs inherent to circular bioeconomic systems.

\subsection{Sensitivity analysis}
We quantify how the optimized objective responds to parameter changes via a local, one-at-a-time (OAT) study around the calibrated baseline (fixed horizon $T=25$). Let $a$ denote a scalar parameter from
\[
\mathcal{P}=\{P_E,\,P_S,\,\alpha,\,\beta,\,\gamma,\,\delta_S,\,\delta_E,\,\rho\}.
\]
Denote by $J_0$ the optimized objective at baseline and by $J(a)$ the optimized objective when only $a$ is perturbed. We use a symmetric finite difference with relative step $\varepsilon=0.10$ (i.e., $\pm10\%$) and report the \emph{relative sensitivity}
\[
S_{\mathrm{rel}}(a)\;:=\;\frac{J\!\left(a(1+\varepsilon)\right)-J\!\left(a(1-\varepsilon)\right)}{2\,\varepsilon\,J_0}\;\cdot\;a.
\]
Hence $S_{\mathrm{rel}}(a)>0$ indicates that a small increase in $a$ raises $J$ (locally), while $S_{\mathrm{rel}}(a)<0$ indicates the opposite; $|S_{\mathrm{rel}}(a)|\approx 1$ means a 1\% change in $a$ produces roughly a 1\% change in $J$. (The cases $\theta{=}0$ vs.\ $\theta{>}0$ and the effect of varying the planning horizon $T$ were analyzed separately in the preceding results sections and are not repeated here.)

\begin{table}[ht]
\centering
\caption{Relative sensitivity of the optimized objective (OAT, $\pm10\%$). Entries are sorted by decreasing $|S_{\mathrm{rel}}|$.}
\label{tab:sensitivity}
\begin{tabular}{@{}lrrrr@{}}
\toprule
\multicolumn{1}{c}{\bf Parameter $a$} &
\multicolumn{1}{c}{\bf Baseline} &
\multicolumn{1}{c}{\bf $J(+10\%)$} &
\multicolumn{1}{c}{\bf $J(-10\%)$} &
\multicolumn{1}{c}{\bf $S_{\mathrm{rel}}(a)$} \\
\midrule
$\rho$       & $0.5000$ & $970.3507$ & $484.6474$ & $1.7621$ \\
$\alpha$     & $0.2500$ & $978.8720$ & $479.7091$ & $0.9055$ \\
$P_S$        & $0.8000$ & $749.5559$ & $629.1333$ & $0.6990$ \\
$\gamma$     & $0.2000$ & $594.1146$ & $804.1123$ & $-0.3048$ \\
$P_E$        & $1.0000$ & $698.2009$ & $680.5154$ & $0.1283$ \\
$\beta$      & $0.3500$ & $698.2009$ & $680.5154$ & $0.0449$ \\
$\delta_S$   & $0.0500$ & $645.1094$ & $736.3796$ & $-0.0331$ \\
$\delta_E$   & $0.0300$ & $689.0797$ & $689.0797$ & $0.0000$ \\
\bottomrule
\end{tabular}
\end{table}

\noindent\textit{Interpretation.} Productivity parameters $\rho$ (soil $\to$ residues) and $\alpha$ (residue return to soil) exert the largest positive influence on $J$, followed by $P_S$. Residue decay $\gamma$ is moderately adverse (negative), indicating higher losses depress performance. Energy-side levers $P_E$ and $\beta$ are positive but smaller at baseline; $\delta_S$ is mildly adverse and $\delta_E$ is locally neutral. All perturbations use identical discretization and solver tolerances to isolate parameter effects (see Section~\ref{sec:numerics}).

\paragraph{Policy implications.}
In our formulation, parameters map cleanly to levers that are widely discussed in the soil–energy literature. The soil valuation term $P_S$ corresponds to instruments that remunerate soil-carbon services (e.g., carbon pricing or result-based agri-environmental payments); raising $P_S$ shifts the optimum toward greater residue return to soil, consistent with integrated assessments showing how incentive design affects the economic efficiency of soil-carbon sequestration. \cite{antle2001econometric} \, In turn, the reinvestment channel (captured in our model by the effectiveness of residue/energy reinjection into soil fertility) has concrete agronomic realizations such as biochar return, which is documented to enhance soil properties and long-term carbon stocks; stronger effectiveness raises the shadow value of soil and favors interior allocations that expand soil capital. \cite{lehmann2015biochar} \, Finally, the energy price $P_E$ summarizes market conditions and support schemes on the bioenergy side; in our numerical results, its variation primarily shifts switching times without altering the bang–interior–bang structure.

\section{Conclusion}\label{sec:conc}
This paper introduced a bioeconomic model for the optimal allocation of agricultural crop residues between energy production and soil fertility enhancement. Using Pontryagin’s Maximum Principle and direct numerical optimization, we derived and analyzed optimal control strategies that reflect the trade-offs between short-term energy revenues and long-term soil regeneration. The model incorporated a novel reinvestment mechanism, where a portion of the accumulated energy benefits is recycled into the system to boost soil productivity, capturing a key feedback in circular bioeconomy systems.
Our simulations revealed rich control structures including switching arcs and interior regimes whose characteristics depend strongly on reinvestment efficiency and planning horizon. In particular, the presence of reinvestment enabled more aggressive early-stage energy extraction without compromising long-term ecological integrity. Comparative analysis with the no-reinvestment scenario underscored the strategic advantage of circular reinforcement and highlighted its potential to improve sustainability outcomes.

These findings support the relevance of optimal control tools for designing residue management strategies in agriculture and bioenergy. Future work will extend the model to include stochastic disturbances (e.g., yield shocks or energy price fluctuations), spatial dynamics, and alternative reinvestment pathways such as composting or infrastructure improvements.


\begin{thebibliography}{99}


\bibitem{alvarez2024}
Alvarez, R.,
A quantitative review of the effects of residue removing on soil organic carbon in croplands,
\textit{Soil Tillage Res.} \textbf{240} (2024) 106098.
DOI: 10.1016/j.still.2024.106098.

\bibitem{anikwe2023circular}
Anikwe, M. A. N., \& Ife, K.,
The role of soil ecosystem services in the circular bioeconomy,
\textit{Front. Soil Sci.} \textbf{3} (2023) 1209100.
DOI: 10.3389/fsoil.2023.1209100.

\bibitem{antle2001econometric}
Antle, J. M., Capalbo, S. M., Elliott, E. T., \& Paustian, K. H.,
Economic analysis of agricultural soil carbon sequestration: An integrated assessment approach,
\textit{J. Agric. Resour. Econ.} \textbf{26} (2001) 344–367.

\bibitem{aseev7}
Aseev, S. M., \& Kryazhimskii, A. V.,
The Pontryagin Maximum Principle and optimal economic growth problems,
\textit{Proc. Steklov Inst. Math.} \textbf{257} (2007) 1–255.
DOI: 10.1134/S0081543807020010.  

\bibitem{beillouin2023}
Beillouin, D., et al.,
A global meta-analysis of soil organic carbon in the world’s agricultural land,
\textit{Nat. Commun.} \textbf{14} (2023) 3630.
DOI: 10.1038/s41467-023-39310-7.

\bibitem{caillau22tur}
Caillau, J.-B., Djema, W., Gouzé, J.-L., Maslovskaya, S., \& Pomet, J.-B.,
Turnpike property in optimal microbial metabolite production,
\textit{J. Optim. Theory Appl.} (2022) 1–33.
DOI: 10.1007/s10957-022-02023-0.

\bibitem{jbcaillau22}
Caillau, J.-B., Ferretti, R., Trélat, E., \& Zidani, H.,
An algorithmic guide for finite-dimensional optimal control problems,
in: \textit{Handb. Numer. Anal.} \textbf{24} (2023) 559–626.
Elsevier/North-Holland.
DOI: 10.1016/bs.hna.2022.11.006.

\bibitem{cherkaoui2022bio}
Cherkaoui Dekkaki, O., El Khattabi, N., \& Raissi, N.,
Bioeconomic modeling of household waste recovery,
\textit{Math. Methods Appl. Sci.} \textbf{45} (2022) 468–482.
DOI: 10.1002/mma.7787.

\bibitem{cherkaoui2023acc}
Cherkaoui Dekkaki, O., \& Djema, W.,
Optimal control of a bioeconomic model applied to the recovery of household waste,
in: \textit{Proc. Amer. Control Conf. (ACC)}, San Diego (2023) 2135–2140.
DOI: 10.23919/ACC55779.2023.10156431.

\bibitem{cherkaoui2024ijdc}
Cherkaoui-Dekkaki, O., Djema, W., Raissi, N., \& El Khattabi, N.,
Optimal control of waste recovery process,
\textit{Int. J. Dyn. Control} \textbf{12} (2024) 4244–4262.
DOI: 10.1007/s40435-024-01484-7.

\bibitem{dekkakiUnderReview}
Cherkaoui Dekkaki, O., Bouhmady, A., \& Benamara, I.,
Multi-objective Optimal Control of a Bioenergy–Soil System under Fertility Constraints,
\textit{NonLinear Science} (2025).


\bibitem{clark2010mathematical}
Clark, C. W., \textit{Mathematical Bioeconomics: The Mathematics of Conservation} (3rd ed.),
Wiley, Hoboken, 2010.

\bibitem{clarke13}
Clarke, F. H., \textit{Functional Analysis, Calculus of Variations and Optimal Control},
Springer, London, 2013.
DOI: 10.1007/978-1-4471-4820-3.

\bibitem{cleanEUislands}
Clean energy for EU islands (European Commission),
Feed-In Premium (Greece)   reference prices for RES/CHP (biomass/biogas),
Available at: \url{https://clean-energy-islands.ec.europa.eu/}
(accessed 8 Oct 2025).


\bibitem{dhanoa2022}
Dhanoa, M. S., Sanderson, R., Cardenas, L. M., Shepherd, A., Chadwick, D. R., Lopez, S., \& France, J.,
Overview and application of the Mitscherlich equation and its extensions to estimate the soil nitrogen pool fraction associated with crop yield and nitrous oxide emission,
\textit{Adv. Agron.} \textbf{174} (2022) 269–295.
DOI: 10.1016/bs.agron.2022.03.005.



\bibitem{dj21turn}
Djema, W., Giraldi, L., Maslovskaya, S., \& Bernard, O.,
Turnpike features in optimal selection of species represented by quota models,
\textit{Automatica} \textbf{132} (2021) 109804.
DOI: 10.1016/j.automatica.2021.109804.

\bibitem{dj22op}
Djema, W., Bayen, T., \& Bernard, O.,
Optimal Darwinian selection of microorganisms with internal storage,
\textit{Processes} \textbf{10} (2022) 461.
DOI: 10.3390/pr10030461.

\bibitem{esma2024}
European Securities and Markets Authority (ESMA),
\textit{Carbon Markets Report 2024},
Available at: \url{https://www.esma.europa.eu/sites/default/files/2024-10/ESMA50-43599798-10379_Carbon_markets_report_2024.pdf}
(accessed 8 Oct 2025).

\bibitem{eurcap2022}
European Union, Council Regulation (EU) 2022/1854   emergency intervention to address high energy prices:
market revenue cap at \texteuro180/MWh for inframarginal generators,
EUR-Lex summary (2022),
Available at: \url{https://eur-lex.europa.eu/}
(accessed 8 Oct 2025).

\bibitem{flichman2011bioeconomic}
Flichman, G., \textit{Bio-Economic Models Applied to Agricultural Systems},
Springer, Dordrecht, 2011.

\bibitem{geshkovski2022}
Geshkovski, B., \& Zuazua, E.,
Turnpike in optimal control of PDEs, ResNets, and beyond,
\textit{Acta Numer.} \textbf{31} (2022) 135–263.
DOI: 10.1017/S0962492922000046.


\bibitem{gordon1954economic}
Gordon, H. S.,
The economic theory of a common-property resource: The fishery,
\textit{J. Polit. Econ.} \textbf{62} (1954) 124–142.
DOI: 10.1086/257497.

\bibitem{rapaport19}
Haddon, A., Ramírez, H., \& Rapaport, A.,
Optimal and sub-optimal feedback controls for biogas production,
\textit{J. Optim. Theory Appl.} \textbf{183} (2019) 642–670.
DOI: 10.1007/s10957-019-01570-3.  

\bibitem{kamien2012dynamic}
Kamien, M. I., \& Schwartz, N. L.,
\textit{Dynamic Optimization: The Calculus of Variations and Optimal Control in Economics and Management},
Dover, Mineola, 2012.

\bibitem{lal2005}
Lal, R.,
World crop residues production and implications of its use as a biofuel,
\textit{Environ. Int.} \textbf{31} (2005) 575–584.
DOI: 10.1016/j.envint.2004.09.005.

\bibitem{lehmann2015biochar}
Lehmann, J., \& Joseph, S. (Eds.),
\textit{Biochar for Environmental Management: Science, Technology and Implementation} (2nd ed.),
Routledge, London/New York, 2015.
DOI: 10.4324/9780203762264.

\bibitem{liu2016optimal}
Liu, Y., \& Liu, M.,
Optimal control of crop irrigation with limited water resources,
\textit{Math. Probl. Eng.} \textbf{2016} (2016) Article ID 5076879, 1–9.
DOI: 10.1155/2016/5076879.

\bibitem{menichetti2019}
Menichetti, L., \AA gren, G. I., Barr\'e, P., Moyano, F., \& K\"atterer, T.,
Generic parameters of first-order kinetics accurately describe soil organic matter decay in bare fallow soils over a wide edaphic and climatic range,
\textit{Sci. Rep.} \textbf{9} (2019) 20319.
DOI: 10.1038/s41598-019-55058-1.

\bibitem{Nunes2024}
Nunes, L. J. R.,
Mathematical modeling approach to the optimization of biomass storage park management,
\textit{Systems} \textbf{12} (2024) 17.
DOI: 10.3390/systems12010017.


\bibitem{perlack2011us}
Perlack, R. D., \& Stokes, B. J.,
\textit{U.S. Billion-Ton Update: Biomass Supply for a Bioenergy and Bioproducts Industry},
Oak Ridge Natl. Lab., 2011.
DOI: 10.2172/1023318.

\bibitem{pont64}
Pontryagin, L. S., \textit{Mathematical Theory of Optimal Processes},
Springer, Berlin/Heidelberg, 1964.


\bibitem{reutersETS2024}
Reuters,
EU carbon market updates (2024–2025): price band roughly \texteuro60–€75/tCO$_2$,
News articles (accessed 8 Oct 2025).


\bibitem{sandoval2024}
Sandoval-Reyes, M., He, R., Semeano, R., \& Ferrão, P.,
Mathematical optimization of waste management systems: Methodological review and perspectives for application,
\textit{Waste Manag.} \textbf{174} (2024) 630–645.
DOI: 10.1016/j.wasman.2023.10.006.


\bibitem{silva2017optimal}
Silva, M., Camargo, M. B. P., \& Schinor, E. H.,
Optimal control in pest management: a case study,
\textit{Ecol. Model.} \textbf{359} (2017) 132–140.

\bibitem{smith2014afolu}
Smith, P., Bustamante, M., Ahammad, H., et al.,
Agriculture, Forestry and Other Land Use (AFOLU),
in: \textit{Climate Change 2014: Mitigation of Climate Change}, Cambridge Univ. Press (2014).
DOI: 10.1017/CBO9781107415416.017.


\bibitem{trelatzuazua2025}
Trélat, E., \& Zuazua, E.,
Turnpike in optimal control and beyond: a survey,
\textit{arXiv preprint} arXiv:2503.20342 (2025).
DOI: 10.48550/arXiv.2503.20342.


\bibitem{tian2024}
Tian, L., Shao, G., Gao, Y., Song, E., \& Lu, J.,
Effects of biochar on soil organic carbon in relation to soil nutrient contents, climate zones and cropping systems: A Chinese meta-analysis,
\textit{Land} \textbf{13} (2024) 1608.
DOI: 10.3390/land13101608.

\bibitem{van1994dynamic}
Van Kooten, G. C., \& Bulte, E. H.,
The economics of nature: managing biological assets,
\textit{Can. J. Agric. Econ.} \textbf{42} (1994) 519–520.

\bibitem{vinter2000}
Vinter, R., \textit{Optimal Control},
Birkhäuser, Boston, 2000.
DOI: 10.1007/978-1-4612-2047-4.

\bibitem{wang2020}
Wang, X., et al.,
Effects of residue returning on soil organic carbon: meta-analysis,
\textit{Agronomy} \textbf{10} (2020) 691.
DOI: 10.3390/agronomy10050691.

\bibitem{yirga2010social}
Yirga, C., \& Hassan, R. M.,
Social costs and incentives for optimal control of soil nutrient depletion in the central highlands of Ethiopia,
\textit{Agric. Syst.} \textbf{103} (2010) 153–160.
DOI: 10.1016/j.agsy.2009.12.002.
\end{thebibliography}
\enddocument